# Low Tensor Train- and Low Multilinear Rank Approximations for De-speckling and Compression of 3D Optical Coherence Tomography Images

Ivica Kopriva, *Senior Member IEEE*, Fei Shi, Mingyig Lai, Marija Štanfel, Haoyu Chen and Xinjian Chen, *Senior Member IEEE*

*Abstract*— Finding optimal representation in tensor decomposition is most often based on prescribed tensor rank or approximation error, both of which are hard to relate to the compression ratio (CR), which is of practical importance. This paper proposes low tensor-train (TT) rank and low multilinear (ML) rank approximations for simultaneous de-speckling and compression of 3D optical coherence tomography (OCT) images for a given CR. Suppression of speckle artifact in OCT imaging is necessary for high-quality quantitative assessment of ocular disorders associated with vision loss. Furthermore, cross-platform software systems for clinical data archiving and/or remote consultation and diagnosis would benefit from compression. To this end, we derive the alternating direction method of multipliers based algorithms for the solution of related de-speckling problems constrained with the low TT- and low ML rank. Rank constraints are implemented through the Schatten-p ($S_p$) norm, $p \in \{0, 1/2, 2/3, 1\}$, of unfolded matrices. We provide the proofs of global convergence towards a stationary point for both algorithms. De-speckled and rank adjusted 3D OCT image tensors are finally approximated through tensor train- and Tucker alternating least squares decompositions. We comparatively validate the low TT- and low ML rank methods on twenty-two 3D OCT images with the JPEG2000 and 3D SPIHT compression methods, as well as with no compression 2D bilateral filtering (BF), 2D median filtering (MF), and enhanced low-rank plus sparse matrix decomposition (ELRpSD) methods. For the CR≤10, the low $S_p$ TT rank method with $p \in \{0, 1/2, 2/3\}$ yields either highest or comparable signal-to-noise ratio (SNR), and comparable or better contrast-to-noise ratio (CNR), mean segmentation errors (SEs) of retina layers and expert-based image quality score (EIQS) than original image and image compression methods. It compares favorably in terms of CNR, fairly in terms of SE and EIQS with the no image compression methods. Thus, for CR≤10 the low $S_{2/3}$ TT rank approximation can be considered a good choice for visual inspection based diagnostics. For 2≤CR≤60, the low $S_1$ ML rank method compares favorably in terms of SE with considered image compression methods and with 2D BF and ELRpSD no compression method. It is slightly inferior to 2D MF no compression method. Thus, for 2≤CR≤60, the low $S_1$ ML rank approximation can be considered a good choice for segmentation based diagnostics either on-site or in the remote mode of operation.

## I. INTRODUCTION

OPTICAL coherence tomography (OCT) is a non-invasive method for measuring reflectance properties of tissue [1]. Thus, quantification of optical properties enables the discrimination of tissues or their pathological states [2] - [6]. Interferometry, the measurement technique upon which OCT is founded [1], [7]-[9], is based on spatial-temporal coherence of the forward and backward scattered signals measured in OCT. However, the same coherence gives rise to speckle, a fundamental property of signals acquired by narrowband detection systems [1]. Speckle reduces contrast and makes boundaries between constitutive tissues more challenging to resolve [7], [1], [10]. That, in return, stands for a significant obstacle in quantitative OCT image analysis [7], [11], [1]. Since OCT itself can be classified as a unique adaptation of electronic speckle pattern interferometry [12], one could argue that complete suppression of speckle is not desirable [1], [13]. That is what makes de-speckling a challenging problem to deal with. De-speckling techniques generally belong to two groups:

This work has been supported by the Croatian Science Foundation Grant IP-2016-06-5235, the European Regional Development Fund under the grant KK.01.1.1.01.0009 (DATACROSS)**,** and the 8[th] Chinese-Croatian Intergovernmental S&T Cooperation Project. The asterisk indicates the corresponding authors. (*Corresponding authors: Ivica Kopriva; Xinjian Chen*).

I. Kopriva is with the Division of Electronics, Ruđer Bošković Institute, Zagreb, Croatia (e-mail: ikopriva@irb.hr).
F. Shi is with the School of Electronics and Information Engineering, Soochow University, Suzhou, China (e-mail:shifei@suda.edu.cn).
M. Lai is with the Shenzen Eye Hospital, Shenzen, China (e-mail: laimydoc@163.com).
M. Štanfel is with the Department of Ophthalmology, University Hospital Center Zagreb, Zagreb, Croatia (e-mail: stanfel.m@gmail.com).
H. Chen is with the Joint Shantou International Eye Center, Shantou University and the Chinese University of Hong Kong, Shantou, China (e-mail: drchenhaoyu@gmail.com).
X. Chen is with the School of Electronics and Information Engineering, Soochow University, Suzhou, China, and the State Key Laboratory of Radiation Medicine and Protection, Soochow University, Suzhou, China (e-mail: xjchen@suda.edu.cn).

physical compounding and digital filtering [8], [14]. Physical compounding methods are hardware-based. They suppress speckle by incoherently summing different realizations of the same OCT image [15]-[18]. These strategies achieve OCT image quality improvement proportional to the square root of the number of realizations. Digital filtering methods aim to suppress speckle through post-processing of OCT image while preserving image resolution, contrast, and edge fidelity [19]-[22]. The digital filtering method can work in the spatial or in the transform domain. Spatial domain methods for speckle suppression range from the traditional ones such as low pass filters, median filters, Wiener filters, and bilinear filters, to the more advanced such as anisotropic filters [23], directional filtering [24] and complex diffusion [25]. Because the OCT image carries on information about the internal microstructure, it is natural to impose a low-rankness constraint on the sought de-speckled image. Low-rankness motivated application of additive matrix decomposition methods to de-speckling of OCT images [22], [26]-[28]. These methods decompose the OCT image on a slice-by-slice basis into the sum of low-rank and sparse matrices [29]. Algorithms such as enhanced low-rank plus sparsity decomposition (ELRpSD) [22] are built upon the improvement of the underlying additive matrix model. Thereby, the low-rank matrix stands for an enhanced image. The exact decomposition is known under the name robust principal component analysis [30] or rank-sparsity decomposition [31].

The de-speckling method proposed herein belongs to the spatial domain group of methods. Its development is motivated by the limitation of additive matrix decomposition methods that presume matrix representation of 3D OCT image. That is, they process a 3D OCT image on a B-scan-by-B-scan basis. However, treating the 3D OCT image as a matrix ignores spatial correlations between the B-scans. The methods proposed herein exploit low-rankness property of the 3D OCT image but, as opposed to matrix representation based methods, de-speckled image is obtained as a result of the minimization of the tensor train (TT) [32], rank and/or multilinear (ML) rank [33]. In addition to the TT representation, a.k.a. matrix product state [34]-[36], the low-rank approximation of 3D OCT image tensor can also be achieved through underlying decompositions such as CANDECOMP/PARAFAC (CP) [37], [38] and Tucker [39]. The tensor rank is not defined uniquely as opposed to matrix rank, but it is decomposition dependent. Computation of the CP rank is proven to be the non-deterministic polynomial-time (NP) hard problem [40]. As discussed in [41]-[43], conceptual drawback of the ML rank, a.k.a. the Tucker rank [41], is that its components are ranks of matrices constructed from and unbalanced tensor matricization scheme (one vs. rest) known as mode-$n$ unfolding, definition 1. Because this scheme results in matrices with a high ratio of maximal and minimal dimension, the rank is effectively equal to the minimal dimension. That can limit the denoising performance of the rank minimization. As opposed to ML rank, the TT rank, defined in [32], consists of ranks of matrices formed by a well-balanced tensor matricization scheme known as mode-(1, 2,..., $k$) or canonical unfolding [33], definition 4. It unfolds the tensor along the permutations of modes. As opposed to the ML rank, the TT rank can capture the global correlation of the tensor entries [41]. Furthermore, mode-(1, 2) unfolding associated with the TT representation of the 3D OCT image tensor contains one matrix less than the mode-$n$ unfolding scheme associated with the ML representation. Consequently, that leads to a computationally more efficient optimization method for low TT rank problem than low ML rank problem, see Figure 3.

As pointed out in [44], the underlying problem in tensor decomposition is to find the optimal representation for the given tensor rank or prescribed approximation error. It is, however, difficult to relate the tensor rank and/or the approximation error to the compression ratio (CR) that is quite often of the practical interest. A trial and error approach is necessary to find approximation error and/or tensor rank that corresponds with the prescribed CR, see Table 1 in [45]. Hence, this paper proposes a two-phase approach to approximation of 3D OCT image tensor by the TT- and/or ML tensor models for the specified CR.

Due to the discontinuous and non-convex nature of the rank function, rank minimization is NP hard problem. Thus, it is often replaced by convex relaxation [30], [46] known as nuclear- or Schatten-1 ($S_1$) norm [47]. $S_1$ norm stands for the sum of singular values and does not represent an accurate measure of rank [48]. A proximity operator associated with the $S_1$ norm over-penalizes large singular values leading to biased results in low-rank constrained optimization problems [49], [50]. To this end, several studies have emphasized the benefit of non-convex penalty functions compared to the nuclear norm for the estimation of the singular values [49], [48], [51]. In addition to that, the selection of the $S_p$ norm, $0 \leq p \leq 1$, is noise-level dependent [52], [53]. $\ell_1$ / $S_1$ norm constrained least-squares problem leads to a more accurate solution when noise is large, while $\ell_p$ / $S_p$, p<1, constraints yield a more accurate solution for a small noise. That, also, is confirmed for the $S_p$ TT- and the $S_p$ ML rank constrained denoising methods proposed herein, see Section V. Discussed arguments stood for motivation to derive alternating direction method of multipliers (ADMM), [54], based algorithms for $S_p$, p$\in$\{0, 1/2, 2/3, 1\}, norm constrained optimization problems associated with minimization of TT- and ML ranks. Selection of these four norms is dictated by the fact that in addition to the soft- and hard thresholding operators associated respectively with the $\ell_1$ norm and $\ell_0$ quasi-norm, the only proximity operators with the closed-form expression for p$\in$(0, 1) exist for p=1/2, [55], and p=2/3, [56].

The main contributions of this paper are as follow:

1) We derived the ADMM-based optimization algorithm for low TT rank approximation of 3D OCT image tensor. The algorithm is applied to matrices constructed from canonical, mode-(1,2,...,$k$), unfolding [33], of 3D OCT image tensor, whereas $S_p$, p$\in$\{0, 1/2, 2/3, 1\}, norms are used to implement low-rank constraints. We provide proof of global convergence towards the stationary point for the derived algorithm.

2) Based on the approximation error between the original and denoised 3D OCT image tensor, the TT-SVD algorithm [32] is used to estimate the core tensors of the TT representation. Information on the obtained TT rank is used to adjust TT rank and generate the final low TT rank approximation of the 3D OCT image with the pre-specified CR, see Section IV.B.

3) We derived the ADMM-based optimization algorithm for low ML rank approximation of 3D OCT image tensor. The algorithm is applied to matrices constructed from mode-$n$ unfolding [33], of 3D OCT image tensor, whereas $S_p$, p∈{0, 1/2, 2/3, 1}, norms are used to implement low-rank constraints. We provide the proof of global convergence towards stationary point for the derived algorithm.

4) ADMM-based low ML rank returns information on ML rank of approximated 3D OCT image tensor. This ML rank is further adjusted to obtain the final low ML rank approximation of the 3D OCT image with the pre-specified CR, see Section IV.C.

5) 3D OCT image tensor denoised by either TT-based or ML-based model can be used for further processing, such as segmentation of the retina layers. It can also be stored in a compressed form. That is useful for cross-platform software systems for clinical data archiving and/or remote consultation and diagnosis. Hence, the proposed low TT rank and/or low ML rank approach to 3D OCT image denoising simultaneously suppress speckle and yield memory efficient low-dimensional representations.

The rest of the paper is organized as follows. Section II presents preliminaries on proximity operators, low-rank approximations, and tensors. Section III presents low TT rank and low ML rank algorithms for denoising of tensors. Section IV details the application of these algorithms to de-speckling and compression of 3D OCT images. Experimental results are presented in Section V. Conclusions are presented in Section VI.

## II. PRELIMINARIES

### A. Notations

We use the following notations throughout the paper. N-dimensional tensor is represented as $\underline{\mathbf{X}} \in \mathbb{R}^{I_1 \times I_2 \times \ldots \times I_N}$ with the elements $x_{i_1 i_2 \ldots i_N}$ where $i_1=1, \ldots, I_1, i_2=1, \ldots, I_2, \ldots, i_N=1, \ldots, I_N$. Each index is called way or mode, and number of levels on one mode is called dimension of that mode, i.e. $\{I_n\}_{n=1}^N$. $\mathbb{R}$ stands for the real manifold. That is standard notation adopted for use in multiway analysis, [57]. Matrices are denoted by uppercase letters in bold font, as an example $\mathbf{X} \in \mathbb{R}^{I_1 \times I_2}$. Vectors are represented with lowercase letters in bold font, as an example $\mathbf{x} \in \mathbb{R}^{I_1}$, while letters in italic font denote scalars, i.e. $x$. The Froebenius norm of $\underline{\mathbf{X}}$ is $\|\underline{\mathbf{X}}\|_F = \sqrt{\sum_{i_1 i_2 \ldots i_N} x_{i_1 i_2 \ldots i_N}^2}$. We also define the $\ell_p$, 0<p≤1, norm of a vector $\ell_p(\mathbf{x}) = \|\mathbf{x}\|_p^p = \sum_{i=1}^N \|x_i\|^p$ and $\ell_0$ quasi-norm of a vector as $\ell_0(\mathbf{x}) = \|\mathbf{x}\|_0 = \#\{x_i \neq 0, i=1,\ldots,N\}$ where # denotes the cardinality function. The $S_p$, 0<p≤1, norms of a matrix $\mathbf{X}$ are defined as corresponding $\ell_p$ norms of the vector of singular values of $\mathbf{X}$, i.e. $S_p(\mathbf{X}) = \|\sigma(\mathbf{X})\|_p^p$, where $\sigma(\mathbf{X})$ stands for the vector of singular values of $\mathbf{X}$. $S_0(\mathbf{X}) = \|\sigma(\mathbf{X})\|_0$.

### B. Proximity operators and low-rank approximation

Let us define the proximity operator of closed proper lower semicontinuous function $g: \mathbb{R} \to \mathbb{R} \cup \infty$ [58], [59]:

$$prox_g^\tau(x) = \arg\min_u \left\{ \frac{1}{2}\|u - x\|^2 + \tau g(u) \right\} \quad (1)$$

For a separable function with the vector argument $\mathbf{x} \in \mathbb{R}^N$ proximal mapping can be computed element wise:

$$prox_g^\tau(\mathbf{x}) = \left[ prox_g^\tau(x_1), \ldots, prox_g^\tau(x_N) \right]^T \quad (2)$$

Low-rank approximation of matrices and tensors is of fundamental importance in many problems in signal processing [60] and machine learning [61]. Owning to the Eckart-Young theorem [62], the low-rank approximation of matrix $\mathbf{X} \in \mathbb{R}^{I_1 \times I_2}$ always exists. It is given in the form of singular-value decomposition (SVD) of $\mathbf{X}=\mathbf{U\Sigma V}$. Low-rank matrix approximation methods approximate $\mathbf{X}$ by minimizing the discrepancy between $\mathbf{X}$ and its model with the appropriate constraints imposed on singular values [49], [63], [64], i.e.

$$G(\mathbf{X}) = \arg\min_\mathbf{A} \left\{ \frac{1}{2}\|\mathbf{X} - \mathbf{A}\|_F^2 + \tau g(\sigma(\mathbf{A})) \right\}$$

where $g(\sigma(\mathbf{A}))$ stands for the rank regularization function. By selecting $g(\sigma(\mathbf{A})) = \|\sigma(\mathbf{A})\|_p^p$ we obtain:

$$G(\mathbf{X}) = \arg\min_\mathbf{A} \left\{ \frac{1}{2}\|\mathbf{X} - \mathbf{A}\|_F^2 + \tau \|\mathbf{A}\|_{S_p} \right\} \quad (3)$$

Note the for the unitary matrices $\mathbf{U}$ and $\mathbf{V}$, $\|\mathbf{UAV}\|_{S_p} = \|\mathbf{A}\|_{S_p}$. Hence, we have [49]:

$$G(\mathbf{X}) = \arg\min_{\mathbf{A}} \left\{ \frac{1}{2} \|\mathbf{X} - \mathbf{A}\|_F^2 + \tau \|\mathbf{A}\|_{S_p} \right\}$$
$$= \mathbf{U} \arg\min_{\mathbf{A}} \left\{ \frac{1}{2} \|\mathbf{\Sigma} - \mathbf{A}\|_F^2 + \tau \|\mathbf{A}\|_{S_p} \right\} \mathbf{V}^T \quad (4).$$

Using the definition and properties of the proximity operators in (1) and (2) we obtain:

$$T_{S_p}(\mathbf{\Sigma}; \tau) = \arg\min_{\mathbf{A}} \left\{ \frac{1}{2} \|\mathbf{\Sigma} - \mathbf{A}\|_F^2 + \tau \|\mathbf{A}\|_{S_p} \right\} \quad (5)$$

where $T_{S_p}(\mathbf{\Sigma}; \tau) = prox_{S_p}^\tau(\mathbf{\Sigma})$ is the $S_p$ norm induced thresholding operator applied entry-wise on $\mathbf{\Sigma}$. Thus, the $S_p$ norm constrained low-rank approximation of $\mathbf{X}$ is given with:

$$\hat{\mathbf{X}} = \mathbf{U} T_{S_p}(\mathbf{\Sigma}; \tau) \mathbf{V}^T \quad (6).$$

For the complete proof of (6) please see the proof of the theorem 2 in [49]. Since $\|\mathbf{X}\|_{S_p} = \|\sigma(\mathbf{X})\|_p^p = \|diag(\mathbf{\Sigma})\|_p^p$, the thresholding operator is $T_{S_p}(\mathbf{\Sigma}; \tau) = T_{\ell_p}(diag(\mathbf{\Sigma}); \tau)$. For $p \in [0, 1]$ analytical expression for $T_{\ell_p}(\circ; \tau)$ exists only for $p \in \{0, 1/2, 2/3, 1\}$. The proximity operators $T_{\ell_0}(\circ, \tau)$ and $T_{\ell_1}(\circ, \tau)$ are, respectively, the well-known hard- and soft-thresholding operators [64]:

$$T_{\ell_0}(x, \tau) = \begin{cases} x, & \text{if } |x| > \sqrt{2\tau} \\ \{0, x\} & \text{if } |x| = \sqrt{2\tau} \\ 0, & \text{if } |x| < \sqrt{2\tau} \end{cases} \quad (7)$$

$$T_{\ell_1}(x, \tau) = sign(x) \max\{|x| - \tau, 0\} \quad (8)$$

The proximity operator $T_{\ell_{1/2}}(\circ, \tau)$ is given with [55]:

$$T_{\ell_{1/2}}(x, \tau) = \begin{cases} \frac{2}{3} x \left(1 + \cos\left(\frac{2\pi}{3} - \frac{2}{3}\alpha(x)\right)\right), & \text{if } |x| > \varphi(\tau) \\ 0, & \text{otherwise} \end{cases} \quad (9)$$

where $\alpha(x) = \arccos(\frac{\tau}{8}(\frac{|x|}{3})^{-\frac{3}{2}})$ and $\varphi(\tau) = \frac{3}{4}(\tau)^{\frac{2}{3}}$. The proximity operator $T_{\ell_{2/3}}(\circ, \tau)$ is given with [56]:

$$T_{\ell_{2/3}}(x, \tau) = \begin{cases} \varpi(x), & \text{if } x > \varphi(\tau) \\ -\varpi(x), & \text{if } x < -\varphi(\tau) \\ 0, & \text{otherwise} \end{cases} \quad (10)$$

with $\varpi(x) = \frac{1}{8}\left(c + \sqrt{\frac{2|x|}{c} - c^2}\right)^3$, $c = \frac{2}{\sqrt{3}} \tau^{\frac{1}{4}} \left(\cosh\left(\frac{\alpha(x)}{3}\right)\right)^{\frac{1}{2}}$, $\alpha(x) = \text{arccosh}\left(\frac{27x^2}{16}\tau^{-\frac{3}{2}}\right)$ and $\varphi(\tau) = \frac{2}{3}(3\tau^3)^{\frac{1}{4}}$.

Usage of $T_{S_0}(\circ, \tau)$, $T_{S_{1/2}}(\circ, \tau)$ and/or $T_{S_{2/3}}(\circ, \tau)$ proximity operators in the low-rank approximation is motivated by the well-known fact that $T_{S_1}(\circ, \tau)$ proximity operator over-penalizes large singular values. That leads to biased results in low-rank constrained optimization problems [49], [50]. That is illustrated in Fig. 1. That in combination with the fact that selection of $\ell_p / S_p$ depends on the noise level motivated us to derive ADMM-based and $S_p$-norm constrained, $p \in \{0, 1/2, 2/3, 1\}$, optimization methods for the low-rank tensor approximation problem.

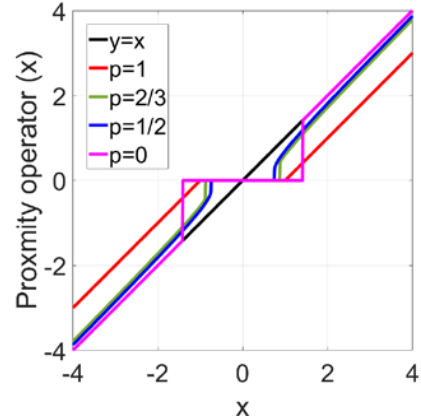

Fig. 1. Proximity operators for threshold value $\tau=1$ associated with the $\ell_p$, $p \in \{0, 1/2, 2/3, 1\}$, norms of x.

### C. Tensors

**Definition 1**. The mode-$n$ matricization, a.k.a. mode-$n$ unfolding, reorders elements of the $N$th order tensor $\underline{\mathbf{X}} \in \mathbb{R}^{I_1 \times I_2 \times \ldots \times I_N}$ for a fixed index $n \in \{1, 2, ..., N\}$ into a matrix:

$$\mathbf{X}_{(n)} \in \mathbb{R}^{I_n \times I_1 I_2 \ldots I_{n-1} I_{n+1} \ldots I_N} \quad (11)$$

Tensor element $(i_1, ..., i_{n-1}, i_n, i_{n+1}, ..., i_N)$ is mapped to the matrix element $(i_n, j)$ such that:

$$j = 1 + \sum_{k=1, k \neq n} (i_k - 1) J_k \text{ with } J_k = \prod_{m=1, m \neq n}^{k-1} I_m \quad (12)$$

The mode-$n$ matricization is implemented in MATLAB through $\mathbf{X}_{(n)} = reshape(\underline{\mathbf{X}}, I_n, I_1 \times ... \times I_{n-1} \times I_{n+1} \times ... \times I_N)$. Conversely, the mode-$n$ matrix can be transformed back to the original tensor by $\underline{\mathbf{X}} = reshape(\mathbf{X}_{(n)}, I_1, I_2, ..., I_N)$.

**Definition 2**. The mode-$n$ product $\underline{\mathbf{C}} = \underline{\mathbf{A}} \times_n \mathbf{B}$ of a tensor $\underline{\mathbf{A}} \in \mathbb{R}^{J_1 \times J_2 \times \ldots \times J_N}$ and a matrix $\mathbf{B} \in \mathbb{R}^{I_n \times J_n}$ is a tensor $\underline{\mathbf{C}} \in \mathbb{R}^{J_1 \times \ldots \times J_{n-1} \times I_n \times J_{n+1} \times \ldots \times J_N}$, with the elements $c_{j_1, j_2 \ldots j_{n-1}, i_n, j_{n+1}, \ldots, j_N} = \sum_{j_n=1}^{J_n} a_{j_1, j_2, \ldots, j_N} b_{i_n, j_n}$.

**Definition 3**. The Tucker decomposition [39] decomposes a $\underline{\mathbf{X}} \in \mathbb{R}^{I_1 \times I_2 \times \ldots \times I_N}$ tensor into a core tensor $\underline{\mathbf{G}} \in \mathbb{R}^{R_1 \times R_2 \times \ldots \times R_N}$ multiplied with factor matrices $\{\mathbf{A}^{(n)} \in \mathbb{R}^{I_n \times R_n}\}_{n=1}^{N}$ on each mode:

$$\underline{\mathbf{X}} \cong \underline{\mathbf{G}} \times_1 \mathbf{A}^{(1)} \times_2 \mathbf{A}^{(2)} \ldots \times_N \mathbf{A}^{(N)} \quad (13)$$

($R_1$, $R_2$, ..., $R_N$) is called the multilinear rank of Tucker decomposition [33].

**Definition 4**. The mode-(1, 2,..., $n$) matricization, a.k.a. mode-$n$ canonical unfolding, reorders elements of the $N$th order tensor $\underline{\mathbf{X}} \in \mathbb{R}^{I_1 \times I_2 \times \ldots \times I_N}$ for a fixed index $n \in \{1, 2, ..., N-1\}$ into a matrix [32], [33]:

$$\mathbf{X}_{[n]} \in \mathbb{R}^{I_1 I_2 \ldots I_n \times I_{n+1} I_{n+2} \ldots I_N} \quad (14)$$

The mode-(1, 2, ..., $n$) matricization is implemented in MATLAB through $\mathbf{X}_{[n]} = reshape(\underline{\mathbf{X}}, I_1 \times \ldots \times I_n, I_{n+1} \times \ldots \times I_N)$. Conversely, the mode-$n$ canonical matrix can be transformed back to tensor by $\underline{\mathbf{X}} = fold_n(\mathbf{X}_{[n]})$. That is implemented in MATLAB through $\underline{\mathbf{X}} = reshape(\mathbf{X}_{[n]}, I_1, I_2, \ldots, I_N)$.

**Definition 5**. Contracted product is defined between the last mode of an $N$-order tensor $\underline{\mathbf{A}} \in \mathbb{R}^{I_1 \times I_2 \times \ldots \times I_N}$ and the first mode of an $K$-order tensor $\underline{\mathbf{B}} \in \mathbb{R}^{J_1 \times J_2 \times \ldots \times J_K}$, where $I_N = J_1$, to yield a tensor $\underline{\mathbf{C}} = \underline{\mathbf{A}} \bullet \underline{\mathbf{B}} \in \mathbb{R}^{I_1 \times I_2 \times \ldots \times I_{N-1} \times J_2 \times \ldots \times J_K}$ with the elements: $c_{i_1 i_2 \ldots i_{N-1} j_2 \ldots j_K} = \sum_{i_N=1}^{I_N} a_{i_1 i_2 \ldots i_{N-1} i_N} b_{i_N j_2 \ldots j_K}$.

**Definition 6**. TT representation of a tensor $\underline{\mathbf{X}} \in \mathbb{R}^{I_1 \times I_2 \times \ldots \times I_N}$ is given in terms of contracted tensor products as:

$$\underline{\mathbf{X}} = \underline{\mathbf{G}}^{(1)} \bullet \underline{\mathbf{G}}^{(2)} \bullet \ldots \bullet \underline{\mathbf{G}}^{(N-1)} \bullet \underline{\mathbf{G}}^{(N)} \quad (15)$$

where third order core tensors, also called the TT-cores, are of the sizes $\{R_{n-1} \times I_n \times R_n\}_{n=1}^{N}$ and it is assumed $R_0 = R_N = 1$. Thus, the first and the last core tensor are matrices: $\mathbf{G}^{(1)} \in \mathbb{R}^{I_1 \times R_1}$ and $\mathbf{G}^{(N)} \in \mathbb{R}^{R_{N-1} \times I_N}$. ($R_1$, $R_2$, ..., $R_{N-1}$) is called tensor train rank.

## III. ALGORITHMS FOR LOW-RANK TENSOR DENOISING

Herein, we present ADMM based algorithms for tensor denoising constrained by the low TT rank and the low ML rank.

### A. Low TT rank denoising of a tensor

The following optimization problem addresses low TT rank denoising of a tensor $\underline{\mathbf{X}} \in \mathbb{R}^{I_1 \times I_2 \times \ldots \times I_N}$:

$$\min_{\mathbf{X}_{[k]}} \sum_{k=1}^{N-1} \alpha_k \, rank(\mathbf{X}_{[k]}) \quad (16)$$

where $\alpha_k$ denotes the weight that the TT rank of the matrix $\mathbf{X}_{[k]}$ contributes to, under the condition $\sum_{k=1}^{N-1} \alpha_k = 1$. As in [41], [42] we select:

$$\alpha_k = \frac{\beta_k}{\sum_{k=1}^{N-1} \beta_k} \quad \text{with} \quad \beta_k = \min\left(\prod_{l=1}^{k} I_l, \prod_{l=k+1}^{N} I_l\right). \quad (17)$$

Since rank minimization is the NP-hard problem, we relax it by the $\|\mathbf{X}_{[k]}\|_{S_p}$ due to the reasons elaborated in Section II.B We select p∈{0, 1/2, 2/3, 1}. By introducing the variable splitting $\mathbf{M}_k = \mathbf{X}_{[k]}$ we obtain the optimization problem:

$$\min_{\mathbf{M}_k} \sum_{k=1}^{N-1} \alpha_k \|\mathbf{M}_k\|_{S_p} \quad \text{such that} \quad \mathbf{M}_k = \mathbf{X}_{[k]} \quad (18)$$

The optimization problem (18) translates into the objective function:

$$\min_{\mathbf{M}_k} \sum_{k=1}^{N-1} \alpha_k \|\mathbf{M}_k\|_{S_p} + \frac{\mu}{2} \|\mathbf{X}_{[k]} - \mathbf{M}_k\|_F^2 \quad (19)$$

where μ>0 is the penalty. The augmented Lagrangian of the objective function (19) is:

$$\mathcal{L}_\mu\left(\{\mathbf{X}_{[k]}\}_{k=1}^{N-1}, \{\mathbf{M}_k\}_{k=1}^{N-1}, \{\mathbf{\Lambda}_k\}_{k=1}^{N-1}\right)$$
$$= \sum_{k=1}^{N-1} \alpha_k \|\mathbf{M}_k\|_{S_p} + \langle \mathbf{\Lambda}_k, \mathbf{X}_{[k]} - \mathbf{M}_k \rangle + \frac{\mu}{2} \|\mathbf{X}_{[k]} - \mathbf{M}_k\|_F^2 \quad (20)$$

where $\{\mathbf{\Lambda}_k\}_{k=1}^{N-1}$ are the Lagrange multipliers.

*Update rule for* $\mathbf{M}_k^{l+1}$. Let $\underline{\mathbf{Z}}^0 = \underline{\mathbf{X}}$. Given the $\mathbf{X}_{[k]}^l$ as a mode-$k$ canonical folding of $\underline{\mathbf{Z}}^l$, $\mathbf{\Lambda}_k^l$ and $\mu^l$ the following problem needs to be solved:

$$\min_{\underline{\mathbf{Z}}^l, \mathbf{M}_k} \left\{ \alpha_k \|\mathbf{M}_k\|_{S_p} + \frac{\mu^l}{2} \left\|\mathbf{X}_{[k]}^l - \mathbf{M}_k + \frac{\mathbf{\Lambda}_k^l}{\mu^l}\right\|_F^2 \right\} \quad (21)$$

Following (5) and (6) we obtain:

$$\mathbf{M}_k^{l+1} = \mathbf{U} T_{S_p}\left(\mathbf{\Sigma}, \frac{\alpha_k}{\mu^l}\right) \mathbf{V}^T \qquad (22)$$

where $\mathbf{U}\mathbf{\Sigma}\mathbf{V}^T$ stands for SVD of $\mathbf{X}_{[k]}^l + \mathbf{\Lambda}_k^l/\mu^l$. Lagrange multipliers are updated according to:

$$\mathbf{\Lambda}_k^{l+1} = \mathbf{\Lambda}_k^l + \mu^l\left(\mathbf{X}_{[k]}^l - \mathbf{M}_k^{l+1}\right). \qquad (23)$$

The penalty parameter is updated according to:

$$\mu^{l+1} = \max\left(\rho\mu^l, \mu_{\max}\right). \qquad (24)$$

Denoised tensor at the iteration $l+1$ is obtained as:

$$\underline{\mathbf{Z}}^{l+1} = \sum_{k=1}^{N-1} \alpha_k \, \text{fold}_k\left(\mathbf{M}_k^{l+1}\right) \qquad (25)$$

The algorithm stops when the convergence criterion based on relative error is satisfied:

$$\frac{\left\|\underline{\mathbf{Z}}^{l+1} - \underline{\mathbf{Z}}^l\right\|_F}{\left\|\underline{\mathbf{X}}\right\|_F} \le \varepsilon_r \qquad (26)$$

or when a maximal number of iteration, *itmax*, is reached. In all the experiments reported in Section V we set $\rho$=1.1, $\varepsilon_r$=0.001 and *itmax*=100. The low TT rank ADMM based algorithm for tensor denoising is summarized in the Algorithm 1.

**Algorithm 1** Low TT rank tensor denoising.

  **Input**: The observed data tensor $\underline{\mathbf{X}}$, penalty parameter $\mu^0$, $\mu_{\max}$, $\rho$, $S_p$ norm for p∈{0, 1/2, 2/3, 1}, targeted value of relative error $\varepsilon_r$ in (26), *itmax*.

  **Parameters**: $\{\alpha_k\}_{k=1}^{N-1}$ based on (17).

1: **Initialization**: $\underline{\mathbf{Z}}^0 = \underline{\mathbf{X}}$, $\{\mathbf{\Lambda}_k^0 = \mathbf{0}\}_{k=1}^{N-1}$, $l$=0, $\varepsilon_r^0 = 1.1\varepsilon_r$.

2: **While** ($\varepsilon_r^l > \varepsilon_r$ or $l<itmax$) **do**:
3:   **for** $k$=1 **to** $N$-1
4:     mode-$k$ canonical unfolding of $\underline{\mathbf{Z}}^l$ to get $\mathbf{X}_{[k]}^l$
5:     $\mathbf{U}\mathbf{\Sigma}\mathbf{V}^T = \mathbf{X}_{[k]}^l + \mathbf{\Lambda}_k^l/\mu^l$
6:     $\mathbf{M}_k^{l+1} = \mathbf{U} T_{S_p}\left(\mathbf{\Sigma}, \frac{\alpha_k}{\mu^l}\right)\mathbf{V}^T$
7:     $\mathbf{\Lambda}_k^{l+1} = \mathbf{\Lambda}_k^l + \mu^l\left(\mathbf{X}_{[k]}^l - \mathbf{M}_k^{l+1}\right)$
8:   **end for**
9:   $\mu^{l+1} = \max(\rho\mu^l, \mu_{\max})$
10:  $\underline{\mathbf{Z}}^{l+1} = \sum_{k=1}^{N-1} \alpha_k \, \text{fold}_k\left(\mathbf{M}_k^{l+1}\right)$
11:  $\varepsilon_r^{l+1} = \left\|\underline{\mathbf{Z}}^{l+1} - \underline{\mathbf{Z}}^l\right\|_F / \left\|\underline{\mathbf{X}}\right\|_F$
12:  $l \leftarrow l+1$
13: **End while**

  **Output**: Denoised data tensor: $\hat{\underline{\mathbf{X}}} \leftarrow \underline{\mathbf{Z}}^l$.

### B. Low ML rank denoising of tensor

The following optimization problems addresses low ML rank denoising of a tensor $\underline{\mathbf{X}} \in \mathbb{R}^{I_1 \times I_2 \times \ldots \times I_N}$:

$$\min_{\mathbf{X}_{(k)}} \sum_{k=1}^{N} \delta_k \, \text{rank}\left(\mathbf{X}_{(k)}\right) \qquad (27)$$

where $\delta_k$ denotes the weight that the ML rank of the matrix $\mathbf{X}_{(k)}$ contributes to, under the condition $\sum_{k=1}^{N} \delta_k = 1$. Herein we select:

$$\delta_k = \frac{\gamma_k}{\sum_{k=1}^{N} \gamma_k} \quad \text{with} \quad \gamma_k = \min\left(I_k, \prod_{l=1, l\neq k}^{N} I_l\right). \qquad (28)$$

We relax $\text{rank}\left(\mathbf{X}_{(k)}\right)$ by the $\left\|\mathbf{X}_{(k)}\right\|_{S_p}$, p∈{0, 1/2, 2/3, 1}. By introducing the variable splitting $\mathbf{M}_k=\mathbf{X}_{(k)}$ we obtain the optimization problem:

$$\min_{\mathbf{M}_k} \sum_{k=1}^{N} \alpha_k \left\|\mathbf{M}_k\right\|_{S_p} \quad \text{such that} \quad \mathbf{M}_k = \mathbf{X}_{(k)} \qquad (29)$$

The optimization sequence related to (29) is analogous to the one for the Algorithm 1, eq. (20) to (26). In the interest of space, we summarize the ADMM based algorithm for low ML rank tensor denoising in the Algorithm 2.

**Algorithm 2** Low ML rank tensor denoising.

  **Input**: The observed data tensor $\underline{\mathbf{X}}$, penalty parameter $\mu^0$, $\mu_{\max}$, $\rho$, $S_p$ norm for p∈{0, 1/2, 2/3, 1}, targeted value of relative error $\varepsilon_r$, *itmax*.

  **Parameters**: $\{\delta_k\}_{k=1}^{N}$ based on (28).

1: **Initialization**: $\underline{\mathbf{Z}}^0 = \underline{\mathbf{X}}$, $\{\mathbf{\Lambda}_k^0 = \mathbf{0}\}_{k=1}^{N}$, $l$=0, , $\varepsilon_r^0 = 1.1\varepsilon_r$.

2: **While** ($\varepsilon_r^l > \varepsilon_r$ or $l<itmax$) **do**:
3:   **for** $k$=1 **to** $N$
4:     mode-$k$ unfolding of $\underline{\mathbf{Z}}^l$ to get $\mathbf{X}_{(k)}^l$
5:     $\mathbf{U}\mathbf{\Sigma}\mathbf{V}^T = \mathbf{X}_{(k)}^l + \mathbf{\Lambda}_k^l/\mu^l$
6:     $\mathbf{M}_k^{l+1} = \mathbf{U} T_{S_p}\left(\mathbf{\Sigma}, \frac{\delta_k}{\mu^l}\right)\mathbf{V}^T$
7:     $\mathbf{\Lambda}_k^{l+1} = \mathbf{\Lambda}_k^l + \mu^l\left(\mathbf{X}_{(k)}^l - \mathbf{M}_k^{l+1}\right)$
8:   **end for**
9:   $\mu^{l+1} = \max(\rho\mu^l, \mu_{\max})$
10:  $\underline{\mathbf{Z}}^{l+1} = \sum_{k=1}^{N} \delta_k \, \text{fold}_k\left(\mathbf{M}_k^{l+1}\right)$

11: $\varepsilon_r^{l+1} = \left\| \underline{\mathbf{Z}}^{l+1} - \underline{\mathbf{Z}}^l \right\|_F / \left\| \underline{\mathbf{X}} \right\|_F$

12: $l \leftarrow l+1$

13: **End while**

**Output**: Denoised data tensor: $\hat{\underline{\mathbf{X}}} \leftarrow \underline{\mathbf{Z}}^l$;

ML rank = $\left( rank\left(\mathbf{M}_1^l\right), ..., rank\left(\mathbf{M}_N^l\right) \right)$.

*C. Convergence Analysis*

We prove the global convergence of Algorithms 1 and 2. The proofs are based on verification of the fulfillment of the five conditions established for convergence of ADMM based minimization of non-convex and non-smooth objective function [67]. As can be seen, Algorithms 1 and 2 are comprised of respectively $N-1$ and $N$ parallel optimization problems. Hence, it suffices to prove convergence for one optimization problem in each case.

**Theorem 1**. Let $D^l = \left\{ \mathbf{M}_k^l, \mathbf{X}_{[k]}^l, \mathbf{\Lambda}_k^l \right\}$ be a sequence generated by Algorithm 1 for problem $k=1,...,N$-1. Then for a sufficiently large $\mu$ Algorithm 1 converges globally, i.e. independently of the initialization point.

**Proof**. We rewrite the problem (19) as:

$$\min_{\mathbf{M}_k, \mathbf{X}_{[k]}} f\left(\mathbf{M}_k\right) + h\left(\mathbf{X}_{[k]}\right) \quad \text{subject to:} \mathbf{A}\mathbf{M}_k = \mathbf{B}\mathbf{X}_{[k]}$$

where $\mathbf{A}=\mathbf{I}$, $\mathbf{B}=\mathbf{I}$, $f\left(\mathbf{M}_k\right) = \left(\alpha_k/\mu\right)\left\|\mathbf{M}_k\right\|_{S_p}$, and $h\left(\mathbf{X}_{[k]}\right) = (1/2)\left\|\mathbf{X}_{[k]} - \mathbf{M}_k\right\|_F^2$. We now show that conditions A1 to A5 that guarantee convergence in the non-convex and possibly non-smooth optimization problem (19) are satisfied [67]. The $S_p$ norm is non-negative lower semi-continuous function bounded from below. Thus, $f$ is also non-negative lower semi-continuous and bounded function. Since $h$ is coercive and $f$ is coercive $f + h$ is also coercive. Thus, assumptions A1 and A4 hold. Since $\mathbf{A}=\mathbf{I}$ and $\mathbf{B}=\mathbf{I}$ assumptions A2 and A3 hold. Since $h(\mathbf{X}_{[k]})$ is Lipschitz differentiable, for sufficiently large $\mu$ sequence $D^l$ is bounded and has at least one limit point. Each limit point is a stationary point of $\mathcal{L}_\mu$ (20), that is $0 \in \partial \mathcal{L}_\mu \left(\mathbf{M}_k^*, \mathbf{\Lambda}_k^*\right)$. Furthermore, since $S_p$ norm is semi-algebraic function, (19) is the Kurdyka-Łojasiewicz function [67], [58], [68]. Thus [67], Algorithm 1 converges globally to the unique stationary point of (20). □

**Theorem 2**. Let $D^l = \left\{ \mathbf{M}_k^l, \mathbf{X}_{(k)}^l, \mathbf{\Lambda}_k^l \right\}$ be a sequence generated by Algorithm 2 for problem $k=1,...,N$. Then, for a sufficiently large $\mu$ Algorithm 2 converges globally.

**Proof**. The only difference between Algorithms 1 and 2 is in the tensor matricization scheme. The reasoning used for the proof of convergence of the Theorem 2 is the same as the reasoning used for the proof of convergence of Theorem 1. □

## IV. DE-SPECKLING AND COMPRESSION IN 3D OCT IMAGE

We now apply the low TT rank and low ML rank denoising algorithms for de-speckling and compression of the 3D OCT image tensor of a retina, i.e. $\underline{\mathbf{X}} \in \mathbb{R}_{0+}^{I_1 \times I_2 \times I_3}$. 3D OCT image tensor is comprised of a set of $I_3$ B-scan images with the size of $I_1 \times I_2$ pixels. The first two modes are for rows and columns, and the third mode is for B-scans. In particular, for 3D OCT images considered in Section V it applies $I_1$=480, $I_2$=512 and $I_3$=64. Before presenting applications of the low TT rank and low ML rank algorithms, we need to introduce the method that aligns B-scans [69].

*A. B-scans alignment method*

The involuntary movement of the eye during OCT acquisition may cause the position of the retina to jump between B-scans and that compromises the B-scan-wise spatial correlations. Therefore, we perform a pre-processing step called B-scan alignment, the same as used in [69]. The position of the upper retinal boundary, detected by a multi-resolution graph search method, is used to estimate the vertical displacement of the retina across B-scans. Specifically, a mean height value is calculated by averaging the height of the leftmost and the rightmost 20% points on the upper retinal boundary. The center part is excluded from the above calculation considering the fact that the concave-shaped fovea appears in the center of some B-scans. Each B-scan is shifted up or down so that the mean height values of the upper peripheral boundary become the same for all B-scans. Experiments showed that this pre-processing step significantly improved the de-speckling and compression performance of the proposed method.

*B. Low TT rank for de-speckling and compression of 3D OCT image tensor*

When applied to the 3D OCT image tensor $\underline{\mathbf{X}}$, the low TT rank algorithm 1 returns denoised tensor $\hat{\underline{\mathbf{X}}}$. Let ε stand for an error between $\underline{\mathbf{X}}$ and its approximation $\hat{\underline{\mathbf{X}}}$:

$$\varepsilon = \left\| \underline{\mathbf{X}} - \hat{\underline{\mathbf{X}}} \right\|_F / \left\| \underline{\mathbf{X}} \right\|_F \quad (30).$$

We now obtain the low TT rank approximation of the 3D OCT image tensor $\underline{\mathbf{X}}$ by:

$$\underline{\mathbf{Y}} = \text{TT-SVD}\left(\underline{\mathbf{X}}, \varepsilon\right) \quad (31).$$

Here TT-SVD stands for the SVD-based TT decomposition algorithm in [32] that for given $\varepsilon$ yields cores of the TT representation (15). Based on definition 6, $\underline{\mathbf{Y}}$ admits TT representation (15) with $N=3$. Hence, for the obtained TT rank $\left(R_1, R_2\right)$, CR of the TT-approximation is given with:

$$CR(R_1, R_2) = \frac{I_1 I_2 I_3}{I_1 R_1 + R_1 I_2 R_2 + R_2 I_3} \quad (32)$$

Since $\hat{\underline{\mathbf{X}}}$ is obtained with the low TT rank denoising Algorithm 1, $CR(R_1, R_2)$ is mostly only close to the pre-specified CR. Thus, we need to correct the TT rank $(R_1, R_2)$. We propose the following correction scheme:

If $CR(\varepsilon) < CR$:

$$R_i = \max(R_1, R_2) \quad (33)$$

else if $CR(\varepsilon) > CR$:

$$R_i = \min(R_1, R_2) \quad (34).$$

If $R_i = R_1$ we correct $R_1$:

$$R_1 = round\left(\frac{I_1 \cdot I_2 \cdot I_3 - CR \cdot R_2 \cdot I_3}{CR \cdot I_1 + CR \cdot I_2 \cdot R_2}\right) \quad (35)$$

$$R_1 = \min(R_1, I_1)$$

else if $R_i = R_2$ we correct $R_2$:

$$R_2 = round\left(\frac{I_1 \cdot I_2 \cdot I_3 - CR \cdot R_1 \cdot I_1}{CR \cdot I_3 + CR \cdot I_2 \cdot R_1}\right) \quad (36).$$

$$R_2 = \min(R_2, I_2)$$

We now obtain the final low TT rank approximation of the 3D OCT image tensor $\underline{\mathbf{X}}$ by:

$$\hat{\underline{\mathbf{X}}} = \text{TT-SVD}\left(\underline{\mathbf{X}}, (R_1, R_2)\right) \quad (37).$$

Now CR of the low TT rank approximation $\hat{\underline{\mathbf{X}}}$ equals the pre-specified value CR. As can be seen from (22), value of the penalty parameter $\mu^l$ influences the threshold value. Thus, for the given $\varepsilon_r$ in Algorithm 1 the initial value $\mu^0$ as well as the maximal value $\mu_{max}$, see (24), play a role in performance of the low TT rank tensor denoising algorithm. It is evident that values of $\mu^0$ and $\mu_{max}$ depend on the choice of $S_p$ norm. Hence, for $\varepsilon_r = 0.001$, $p \in \{0, 1/2, 2/3, 1\}$ we used twenty two 3D OCT images to tune values of $\mu^0$ and $\mu_{max}$ so that CR (32) ranges approximately from 1 to 100 in twelve steps. Fig 2. shows values of $\mu^0$ vs. mean values of CR. For each selection of the $S_p$ norm we obtain $\mu^0$, as a function of CR, by the shape-preserving cubic Hermite interpolation using MATLAB function `pchip`. Thus, for the selected value of CR we run the low TT rank tensor denoising Algorithm 1 and TT-SVD algorithms (31) and (37) to obtain approximated tensor $\hat{\underline{\mathbf{X}}}$ with the pre-specified CR. Fig. 3 illustrates computation times of the sequence comprised of Algorithm 1, TT-SVD (31) and TT-SVD (37), and the sequence comprised of Algorithm 2 and Tcuker_ALS (44), see Section IV.C, averaged over twenty-two 3D OCT images for selected values of CR. It is seen that low TT rank based approximation is computationally always more efficient than low ML rank approximation described in Section IV.C.

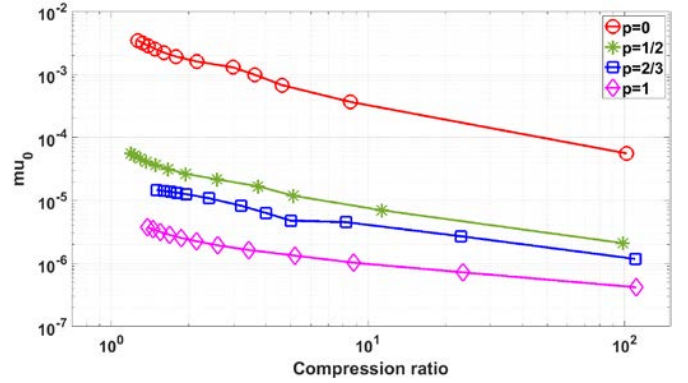

Fig. 2. Initial values of the penalty parameter $\mu^0$ in the low TT rank denoising Algorithm 1, associated with the choice of $S_p$, $p \in \{0, 1/2, 2/3, 1\}$, norm as a function of the CR.

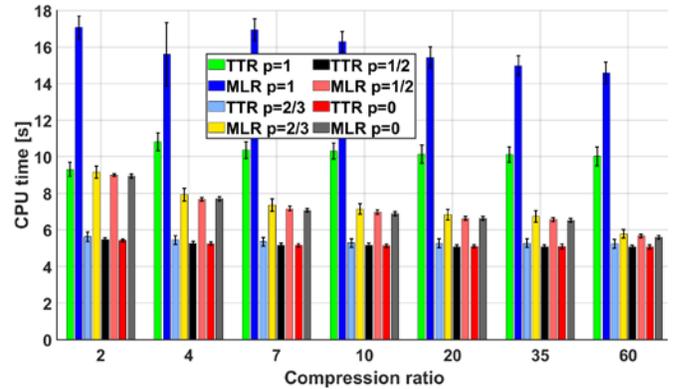

Fig. 3. Computation times (mean ± standard deviation) of the low TT rank and low ML rank tensor denoising algorithms averaged over twenty-two 3D OCT images vs. CR, $S_p$, $p \in \{0, 1/2, 2/3, 1\}$. Stopping criterion: $\varepsilon_r \leq 0.001$ for the low TD rank algorithm and $\varepsilon_r \leq 0.003$ for the low ML rank algorithm.

We summarize the low TT rank approach to de-speckling and compression of 3D OCT images in the Algorithm 3.

**Algorithm 3** Low TT rank approach to de-speckling and compression of 3D OCT images.
**Input**: The observed data tensor $\underline{\mathbf{X}}$, CR, $\rho$, $S_p$ norm for $p \in \{0, 1/2, 2/3, 1\}$, targeted value of relative error $\varepsilon_r = 0.001$ in (26), $itmax = 100$, CR.
**Parameters**: $\{\alpha_k\}_{k=1}^2$ based on (17).
1: Interpolate $\mu^0$ and $\mu_{max}$ based on CR.
2: $\hat{\underline{\mathbf{X}}}$ = Algorithm 1 ($\underline{\mathbf{X}}$, $\mu^0$, $\mu_{max}$, $\rho$, $S_p$, $\varepsilon_r$, $\{\alpha_k\}_{k=1}^2$, $itmax$)
3: Approximation error $\varepsilon = \|\underline{\mathbf{X}} - \hat{\underline{\mathbf{X}}}\|_F / \|\underline{\mathbf{X}}\|_F$ in (30).

4: $\underline{\mathbf{Y}} = \text{TT-SVD}(\underline{\mathbf{X}}, \varepsilon)$ in (31).
5: $(R_1, R_2)$ = TT rank of $\underline{\mathbf{Y}}$.
6: Corrected $(R_1, R_2)$ based on (32) to (36).
7: The final approximation with pre-specified CR:
   $\hat{\underline{\mathbf{X}}} = \text{TT-SVD}(\underline{\mathbf{X}}, (R_1, R_2))$ in (37).

*C. Low ML rank for de-speckling and compression of the 3D OCT image tensor*

Low ML rank algorithm, when applied to the 3D OCT image tensor $\underline{\mathbf{X}}$, returns denoised tensor $\hat{\underline{\mathbf{X}}}$. In addition to that, the low ML rank denoising method returns estimated ML rank=($R_1$, $R_2$, $R_3$). By following definition 3, we can estimate CR for obtained ML rank as:

$$CR(R_1, R_2, R_3) = \frac{I_1 I_2 I_3}{R_1 R_2 R_3 + I_1 R_1 + I_2 R_2 + I_3 R_3} \quad (38).$$

Since $\hat{\underline{\mathbf{X}}}$ is obtained with the low ML rank denoising Algorithm 2, $CR(R_1, R_2, R_3)$ is mostly only close to the pre-specified CR. Thus, we need to correct the ML rank=($R_1$, $R_2$, $R_3$):

If $CR(R_1, R_2, R_3) < CR$:

$$R_i = \max(R_1, R_2, R_3) \quad (39)$$

else if $CR(R_1, R_2, R_3) > CR$:

$$R_i = \min(R_1, R_2, R_3) \quad (40).$$

If $R_i = R_1$ we correct $R_1$:

$$R_1 = \text{round}\left(\frac{I_1 \cdot I_2 \cdot I_3 - CR \cdot R_2 \cdot I_2 - CR \cdot I_3 \cdot R_3}{CR \cdot R_2 \cdot R_3 + CR \cdot I_1}\right) \quad (41)$$

$$R_1 = \min(R_1, I_1)$$

else if $R_i = R_2$ we correct $R_2$:

$$R_2 = \text{round}\left(\frac{I_1 \cdot I_2 \cdot I_3 - CR \cdot I_1 \cdot R_1 - CR \cdot I_3 \cdot R_3}{CR \cdot R_1 \cdot R_3 + CR \cdot I_2}\right) \quad (42)$$

$$R_2 = \min(R_2, I_2)$$

else if $R_i = R_3$ we correct $R_3$:

$$R_3 = \text{round}\left(\frac{I_1 \cdot I_2 \cdot I_3 - CR \cdot I_1 \cdot R_1 - CR \cdot I_2 \cdot R_2}{CR \cdot R_1 \cdot R_2 + CR \cdot I_3}\right) \quad (43).$$

$$R_3 = \min(R_3, I_3)$$

We now obtain the final low ML rank approximation of the 3D OCT image tensor $\underline{\mathbf{X}}$ using the Tucker ALS approach [40], [70], [71]:

$$\hat{\underline{\mathbf{X}}} = \text{Tucker\_ALS}(\underline{\mathbf{X}}, (R_1, R_2, R_3)) \quad (44).$$

Now CR of the low ML rank approximation $\hat{\underline{\mathbf{X}}}$ equals the pre-specified value of the CR. By using the same procedure as the one explained in Section IV.B we tuned values of $\mu^0$ and $\mu_{\max}$ in the low ML rank tensor denoising algorithm such that CR (38) ranges approximately from 1 to 100 in twelve steps. Fig 4. shows values of $\mu^0$ vs. mean values of the CR. For a CR's query value, the value of $\mu^0$ is obtained by the shape-preserving cubic Hermite interpolation using MATLAB function `pchip`.

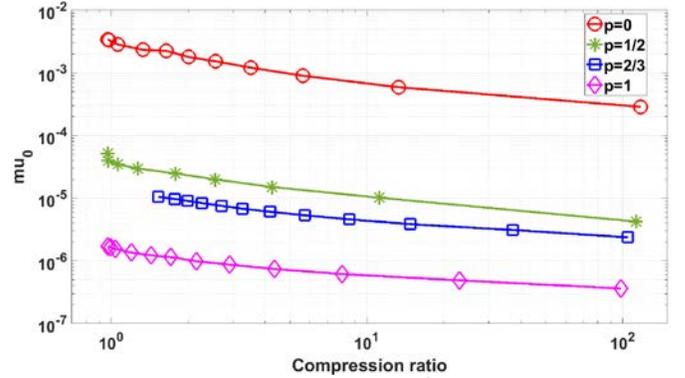

Fig. 4. Initial values of the penalty parameter $\mu^0$ in the low ML rank tensor denoising Algorithm 2, associated with the $S_p$, $p \in \{0, 1/2, 2/3, 1\}$, norm a, as a function of the CR.

For each selection of the $S_p$ norm we obtain $\mu^0$ as a function of CR. Our initial attempt was to stop the low ML rank algorithm when the convergence criterion $\varepsilon_r = 0.001$ was fulfilled. That, however, was not achievable. Hence, we had to increase the value of the convergence criterion to $\varepsilon_r = 0.003$. As can be seen in Fig. 3, low ML rank approach to de-speckling and compression is computationally more demanding than low TT rank approach. That is a consequence of different matricization schemes because of which low ML rank algorithm needs to solve one optimization problem more than the low TT rank algorithm. We summarize the low ML rank method in Algorithm 4.

**Algorithm 4** Low ML rank approach to de-speckling and compression of 3D OCT image.
**Input**: The observed data tensor $\underline{\mathbf{X}}$, CR, $\rho$, $S_p$ norm for $p \in \{0, 1/2, 2/3, 1\}$, targeted value of relative error $\varepsilon_r = 0.003$ in (26), *itmax*=100, CR.

**Parameters**: $\{\delta_k\}_{k=1}^3$ based on (28).

1: Interpolate $\mu^0$ and $\mu_{\max}$ based on CR.

2: $\hat{\underline{\mathbf{X}}}$ = Algorithm 2 ($\underline{\mathbf{X}}$, $\mu^0$, $\mu_{\max}$, $\rho$, $S_p$, $\varepsilon_r$, $\{\delta_k\}_{k=1}^3$, *itmax*)

3: $(R_1, R_2, R_3)$ is estimated ML rank of $\hat{\underline{\mathbf{X}}}$.
4: Corrected $(R_1, R_2, R_3)$ based on (38) to (43).
5: The final approximation with the pre-specified CR:
$\hat{\underline{\mathbf{X}}} = \text{Tucker\_ALS}(\underline{\mathbf{X}}, [R_1, R_2, R_3])$ in (44).

## V. EXPERIMENTAL RESULTS

### A. Data

De-speckling and image compression algorithms were tested comparatively on twenty-two 3D macular-centered OCT images of normal eyes acquired with the Topcon 3D OCT-1000 scanner. Each 3D OCT image was comprised of 64 2D scans with a size of 480×512 pixels. These images were used previously for the study for optical intensity analysis in [72], where they were segmented into retina layers. They are publicly available at [73]. The images were analyzed with the software written in the MATLAB® (the MathWorks Inc., Natick, MA) script language on PC with Intel i7 CPU with the clock speed of 2.2 GHz and 16GB of RAM. MATLAB implementation of the TT-SVD algorithm was downloaded from [74].

### B. Other methods for comparison

We compared the low TT rank and low ML rank algorithms for de-speckling and compression of 3D OCT images with the JPEG2000 standard for image compression, the 3D SPIHT image compression method [75], and no compression de-speckling methods: 2D bilateral filtering (BF), 2D median filtering (MF) and ELRpSD [22]. For BF, the parameters are chosen so that the mean contrast-to-noise ratio (CNR) is comparable to the low TT rank tensor denoising algorithm with CR=10. We used the *imwrite* function in MATLAB to perform JPEG2000 compression. The 3D SPIHT method for 2 bytes/pixel images was downloaded as the executable program from [76]. We set the CRs of these two methods to be the same as the CRs selected for the low TT rank and low ML rank methods. The ELRpSD method was compared extensively in [22] with additive low-rank plus sparse matrix decomposition methods for de-speckling such as [77], [78]. The ELRpSD algorithm yielded the best quality in terms of the performance based on a combination of the CNR and signal-to-noise ratio (SNR) measures. Thus, it is fair to use the ELRpSD as a representative for low-rank sparse additive matrix decomposition-based de-speckling methods.

### C. Performance measures

In order to quantify the performance of algorithms for de-speckling and compression appropriate measures have to be defined. First, we segmented retina layers from each enhanced 3D OCT image as well as from original images. Since de-speckling is usually applied as a pre-processing step of automatic OCT image analysis, such as segmentation, we would like to evaluate how it affects the subsequent segmentation result. Thus, we compared corresponding segmented layers with the ones obtained by manual segmentation in [72]. Segmentation error (SE) is calculated as:

$$SE = \frac{1}{I_2 |S| L} \sum_{i_2=1}^{I_2} \sum_{i_3 \in S} \sum_{l=1}^{L} \frac{|AS^l(i_2, i_3) - MS^l(i_2, i_3)|}{T(i_2, i_3)} \quad (45)$$

where $AS^l(i_2, i_3)$ and $MS^l(i_2, i_3)$ denote the vertical position of the $l$th segmenting surface obtained by the automatic method and the manual segmentation, respectively. $T(i_1, i_2)$ stands for the total thickness of the retina at position $(i_2, i_3)$. Since manual segmentation of the retina layers is tedious, we used $L=8$ surfaces and a subset of B-scans with indices belonging to $S = \{i \times 6\}_{i=1}^{10}$ for the SE estimation [72]. Thus, the SE represents the mean absolute vertical difference between the automatic and manual segmentation relative to the total thickness of a retina. Furthermore, in the case of OCT image de-speckling, the most commonly used figure of merit is CNR [8], [7]. It corresponds to the inverse of the speckle fluctuation and it is defined as: $CNR = \mu_l(\underline{\mathbf{X}}) / \sigma_l(\underline{\mathbf{X}})$ where $\mu_l(\underline{\mathbf{X}})$ and $\sigma_l(\underline{\mathbf{X}})$ respectively correspond to the mean and standard deviation in some selected homogeneous part of the image $\underline{\mathbf{X}}$. Experimental results reported in Section V.D were estimated in the region that corresponds with the topmost layer in the OCT image of a retina, indicated in Fig. 15 by an arrow [5], [72]. Since the goal of post-processing algorithms is not only to suppress speckle, but also to preserve image resolution [8], experienced ophthalmologist estimated quality of original and de-speckled OCT images. To this end, mean and standard deviation of the expert image quality score (EIQS) were estimated from the three B-scans of the three 3D OCT images de-speckled by each considered method as well as the original images. The EIQS ranges from 1 to 5, where 1 means the poorest quality and 5 the highest quality.

### D. Results of comparative performance analysis

Fig. 5 shows the mean values of the SE vs. the CR for twenty-two 3D OCT images enhanced with the low TT rank and low ML rank methods constrained with four $S_p$ norms. Fig. 6 shows mean values of the SE vs. the CR for the low $S_{2/3}$ TT rank and low $S_1$ ML rank algorithms, JPEG2000 and 3D SPIHT algorithms. Corresponding results analogous to Fig. 5 are for CNR, SNR and EIQS shown in Figs. 7, 8 and 9. Herein, it is necessary to comment the Figs. 5, 7 and 8 where low $S_p$ TT and ML rank methods achieved the same performance for $p \in \{0, 1/2, 2/3\}$. Evidently, after Algorithm 1 and Algorithm 2, the rank adjustment procedures (32)-(36) and (38)-(43) converged towards very close or same TT- and ML ranks for $p \in \{0, 1/2, 2/3\}$. Thus, the TT-SVD algorithm in (37) and Tucker\_ALS algorithm in (44) yielded, from the SE, CNR and SNR viewpoints, the same or highly similar results. Results analogous to Fig. 6 are for the EIQS vs. show in Fig. 10. Figs. 11 to 14 show respectively the mean values (± standard deviation) of the SE, CNR, SNR and EIQS for methods discussed in Section V.B. For the CR≤10, the low $S_p$ TT rank model with p∈{0, 1/2, 2/3} yields the highest or comparable SNR, and comparable or better CNR, SE and EIQS values than the original image and JPEG2000 and 3D SPIHT image compression methods. It compares favorably in terms of CNR, fairly in terms of SE and EIQS with the no image compression methods 2D BF, 2D MF and 2D ELRpSD. Hence, for the CR≤10 the low $S_{2/3}$ TT rank method can be used for visual assessment of ocular disorders and/or segmentation-based diagnostic either on-site or in remote consultation. For 2≤CR≤60 the low $S_1$ ML rank model compares favorably in terms of the SE with all compression methods and 2D BF and ELRpSD no compression methods and it is slightly inferior to 2D MF no compression method, see Figs. 6 and 11. For 10≤CR≤60, it is inferior to JPEG and 3D SPIHT methods in terms of EIQS, see Fig. 10. Hence, the low $S_1$ ML rank method could be used in segmentation-based diagnostic scenarios for 2≤CR≤60 either in-site or in the remote mode of operation. Discussed results agree with [52], [53], where it was shown that $\ell_1 / S_1$ norm constrained least-squares problem yields more accurate solution when noise is large. In our application, that coincides with the by large quantization error caused by high CR. For the purpose of illustration, Fig. 15 shows enhancement results on one exemplary B-scan from one original 3D OCT image and corresponding 3D OCT images enhanced with: the low $S_{2/3}$ TT rank for CR=7, the low $S_1$ ML rank for CR=60, JPEG2000 for CR=7 and CR=60, 3D SPIHT for CR=10 and CR=60, the ELRpSD, 2D BF and 2D MF algorithms. Computation times of selected 3D OCT image enhancement methods averaged over twenty-two images are given in Table I. With the exception of 3D SPIHT algorithm, that was implemented in C language, all algorithms were implemented in MATLAB. Low $S_{2/3}$ TT rank and low $S_1$ ML rank methods have moderate computational complexity that is justified by performance achieved in comparison with other methods under consideration.

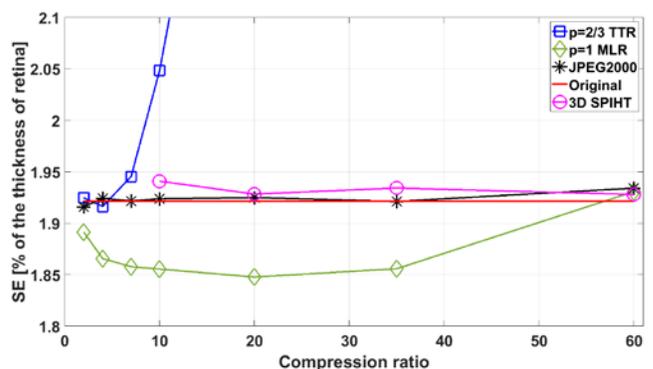

Fig. 6. Best seen in color. Mean values of SE (45) vs. selected values of CR for the $S_{2/3}$ TD-TTR and $S_1$ TD-MLR algorithms, JPEG2000 and 3D SPIHT algorithms.

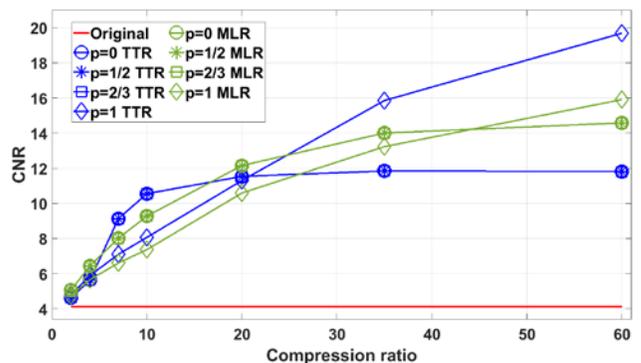

Fig. 7. Best seen in color. Mean values of CNR vs. selected values of CR for the TD-TTR and TD-MLR algorithms for $S_p$, p∈{0, 1/2, 2/3, 1} and 3D OCT image tensors.

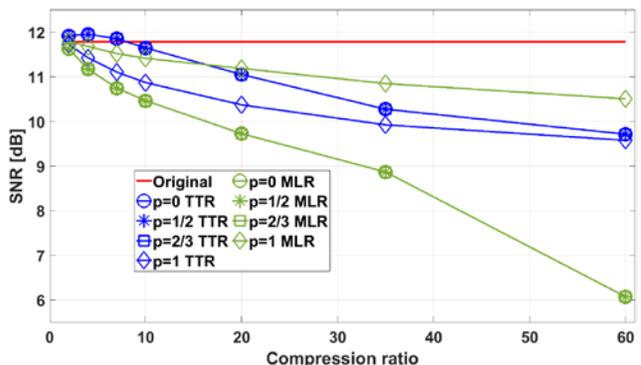

Fig. 8. Best seen in color. Mean values of SNR vs. selected values of CR for the TD-TTR and TD-MLR algorithms for $S_p$, p∈{0, 1/2, 2/3, 1} and 3D OCT image tensors.

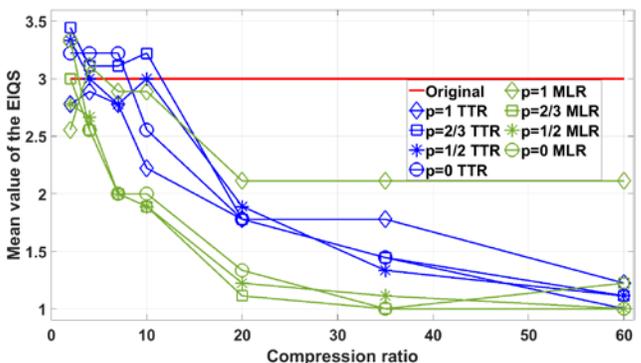

Fig. 9. Best seen in color. Mean values of the EIQS vs. selected values of CR for the TD-TTR and TD-MLR algorithms for $S_p$, p∈{0, 1/2, 2/3, 1} and 3D OCT image tensors.

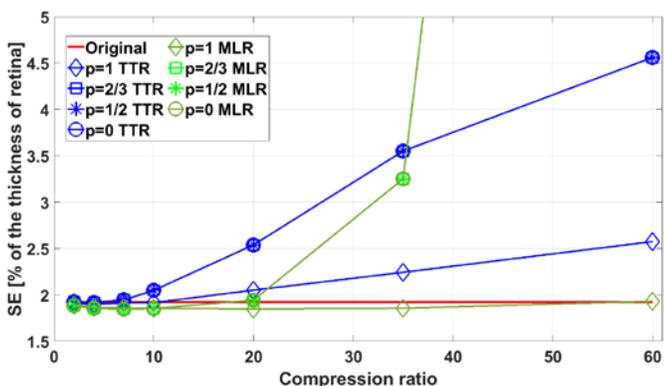

Fig. 5. Best seen in color. Mean values of SE (45) vs. selected values of CR for the TD-TTR and TD-MLR algorithms for $S_p$, p∈{0, 1/2, 2/3, 1} 3D OCT image tensors.

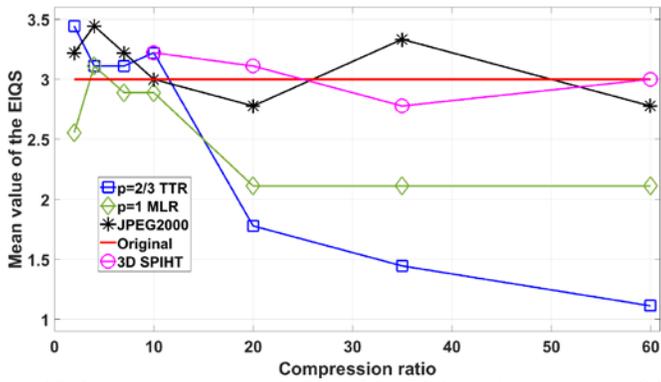

Fig. 10. Best seen in color. Mean values of the EIQS vs. selected values of CR for the TD-TTR and TD-MLR algorithms for $S_p$, p∈{0, 1/2, 2/3, 1} and 3D OCT image tensors.

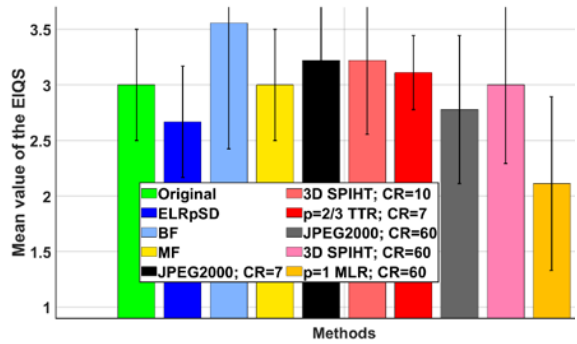

Fig. 14. Best seen in color. EIQS (mean ± standard deviation).

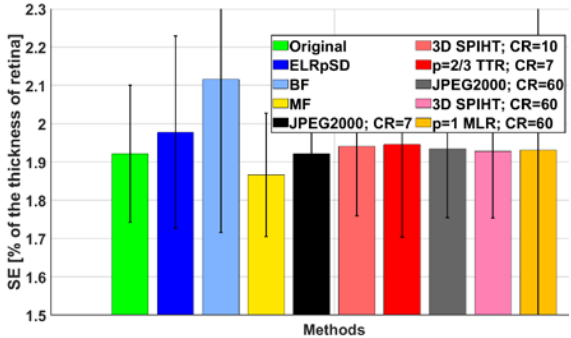

Fig. 11. Best seen in color. Relative SE (mean ± standard deviation).

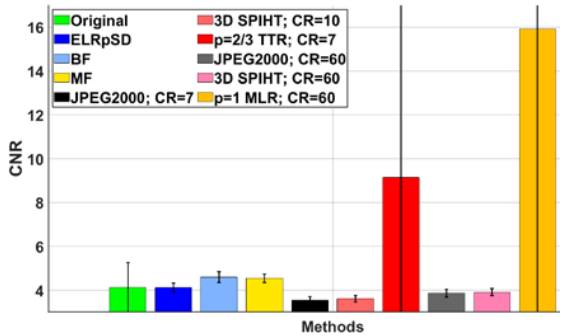

Fig. 12. Best seen in color. CNR (mean ± standard deviation).

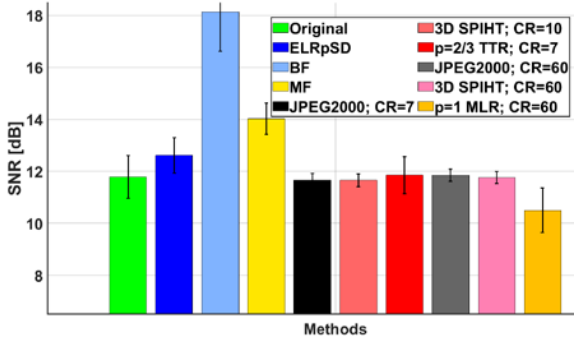

Fig. 13. Best seen in color. SNR (mean ± standard deviation).

TABLE I
COMPUTATION TIMES (MEAN ± STANDARD DEVIATION)

| Method | Computation time [s] |
| --- | --- |
| Low $S_{2/3}$ TT rank (CR=7) | 5.31 ± 0.12 |
| Low $S_1$ ML rank (CR=60) | 14.58 ± 0.60 |
| 3D SPIHT (CR=60) | 1.23 ± 0.02 (C language) |
| JPEG 2000 (CR=60) | 2.64 ± 0.44 |
| 2D BF | 400.32 ± 2.35 |
| 2D MF | 0.426 ± 0.014 |
| ELRpSD | 38.19 ± 0.230 |

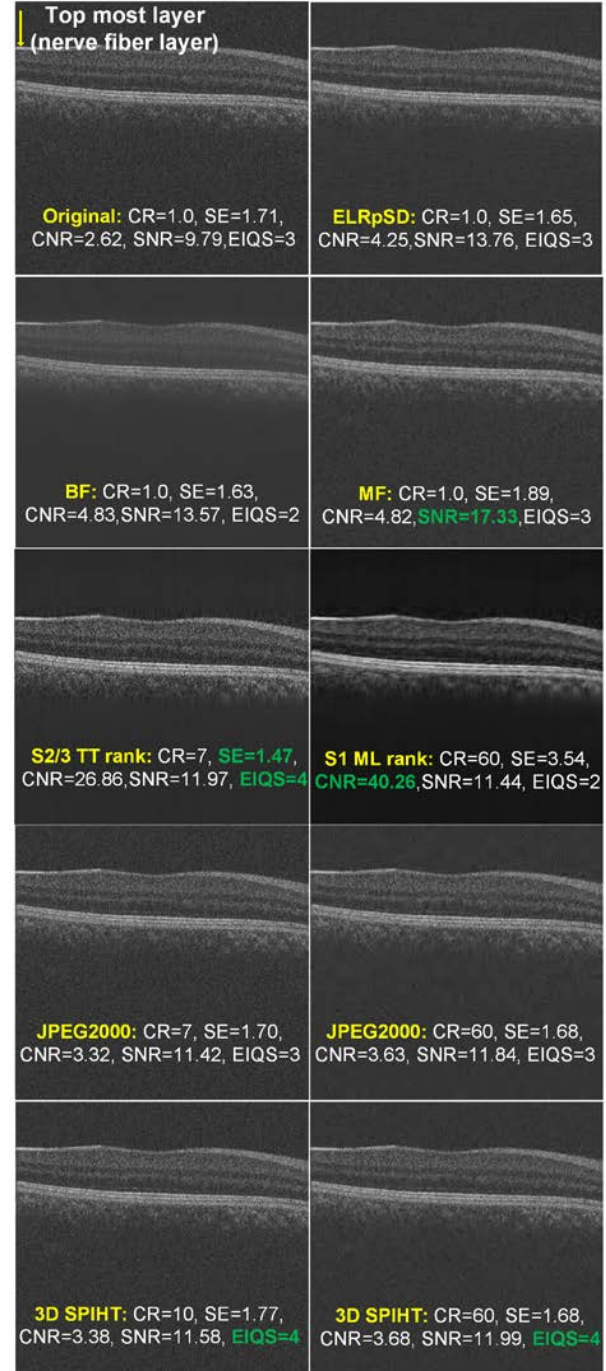

Fig. 15. Best seen in color. 3D OCT image B-scan 27. The best value of each measure is in green bold. For visual comparison, B-scans were mapped to [0, 1] interval with the MATLAB `mat2gray` command from the interval corresponding to minimal and maximal values of the each specific case.

## VI. CONCLUSIONS

Simultaneous de-speckling and compression of 3D OCT images is vital for quantitative assessment of ocular disorders associated with vision loss, clinical data archiving and/or remote consultation and diagnosis. To take into account the tensorial nature of 3D OCT images, this paper proposed low TT rank and low ML rank approximations of 3D OCT image, whereas $S_p$, $p\in\{0, 1/2, 2/3, 1\}$, norm was used as a low rank constraint. To this end, ADMM-based algorithms were derived to solve related optimization problems. In particular, low TT rank and low ML rank approximations of 3D OCT image tensor are obtained for a prescribed CR. This stands for a practically significant contribution because finding optimal representation in tensor decomposition is most often based on prescribed tensor rank or approximation error. Both of them are hard to relate to the CR, which is of practical importance. For CR$\leq$10, comparative performance analysis on twenty-two 3D OCT images of a retina showed that low $S_p$, $p\in\{0, 1/2, 2/3\}$, TT rank algorithms yield 3D OCT images of competitive quality in terms of SNR, CNR, and EIQS. That qualifies them for use in scenarios when the visual assessment of ocular disorders is required either on-site or in remote consultation. On the other side, for 2$\leq$CR$\leq$60 the low $S_1$ ML rank method compares favorably in terms of the SE of the retina layers among compared compression and 2D BF and ELRpSD no compression methods. It is slightly inferior to the 2D MF no compression methods. Hence, it could be useful in machine-learning based diagnostic scenarios, either on-site or in the remote mode of operation.